\documentclass[12pt]{amsart}

\usepackage{latexsym}
\usepackage{amsthm,amsfonts,amssymb,amscd,mathrsfs}
\usepackage[leqno]{amsmath}
\usepackage{xspace}
\usepackage{cite}
\usepackage{fancyhdr}

\newcommand{\numberseries}{\mdseries}	

\newlength{\thmtopspace}		
\newlength{\thmbotspace}		
\newlength{\thmheadspace}		
\newlength{\thmindent}			

\newtheoremstyle{bfupright head,slanted body}
		{\thmtopspace}{\thmbotspace}
		{\slshape}{\thmindent}{\bfseries}{.}{\thmheadspace}
		{{\numberseries \thmnumber{(#2) }}\thmnote{#3}}

\newtheoremstyle{bfupright head,upright body}
		{\thmtopspace}{\thmbotspace}
		{\upshape}{\thmindent}{\bfseries}{.}{\thmheadspace}
		{{\numberseries \thmnumber{(#2) }}\thmnote{#3}}

\newtheoremstyle{bfit head,upright body}
		{\thmtopspace}{\thmbotspace}
		{\upshape}{\thmindent}{\upshape}{.}{\thmheadspace}
		{{\numberseries\thmnumber{(#2) }}
		{\bfseries\itshape\thmnote{\negthickspace#3}}}

\newtheoremstyle{it head,upright body}
		{\thmtopspace}{\thmbotspace}
		{\upshape}{\thmindent}{\upshape}{.}{\thmheadspace}
		{{\numberseries\thmnumber{(#2) }}
		{\itshape\thmnote{\negthickspace#3}}}

\newtheoremstyle{fixed bf head,slanted body}
		{\thmtopspace}{\thmbotspace}{\slshape}
		{\thmindent}{\bfseries}{.}{\thmheadspace}
		{{\numberseries \thmnumber{(#2) }}\thmname{#1}\thmnote{ (#3)}}

\newtheoremstyle{fixed bf head,upright body}
		{\thmtopspace}{\thmbotspace}{\upshape}
		{\thmindent}{\bfseries}{.}{\thmheadspace}
		{{\numberseries \thmnumber{(#2) }}\thmname{#1}\thmnote{ (#3)}}

\newtheoremstyle{indented paragraph}
		{\thmtopspace}{\thmbotspace}
		{\upshape}{\thmindent}{\upshape}{}{0pt}
		{\thmnote{#3 }}

\theoremstyle{bfupright head,slanted body}
\newtheorem{res}{}[section]		\newtheorem*{res*}{}

\theoremstyle{bfit head,upright body}
			\newtheorem*{com*}{}

\theoremstyle{bfupright head,upright body}
\newtheorem{bfhpg}[res]{}		\newtheorem*{bfhpg*}{}

\theoremstyle{it head,upright body}
		\newtheorem*{ithpg*}{}

\theoremstyle{fixed bf head,slanted body}
\newtheorem{thm}[res]{Theorem}		\newtheorem*{thm*}{Theorem}
\newtheorem{prp}[res]{Proposition}	\newtheorem*{prp*}{Proposition}
	\newtheorem*{cor*}{Corollary}
\newtheorem{lem}[res]{Lemma}		\newtheorem*{lem*}{Lemma}

\theoremstyle{fixed bf head,upright body}
\newtheorem{dfn}[res]{Definition}	\newtheorem*{dfn*}{Definition}
\newtheorem{obs}[res]{Observation}	\newtheorem*{obs*}{Observation}
\newtheorem{rmk}[res]{Remark}		\newtheorem*{rmk*}{Remark}
\newtheorem{exa}[res]{Example}		\newtheorem*{exa*}{Example}
		\newtheorem*{exe*}{Exercise}
		\newtheorem{stp*}{Setup}
	\newtheorem*{dfns*}{Definitions}
	\newtheorem*{obss*}{Observations}
		\newtheorem*{rmks*}{Remarks}
	\newtheorem*{exas*}{Examples}

\theoremstyle{indented paragraph}

\newlength{\thmlistleft}	
\newlength{\thmlistright}	
\newlength{\thmlistpartopsep}	
\newlength{\thmlisttopsep}	
\newlength{\thmlistparsep}	
\newlength{\thmlistitemsep}	

\newcounter{eqc} 
\newenvironment{eqc}{\begin{list}{\upshape (\textit{\roman{eqc}})}%
		    {\usecounter{eqc}%
		        \setlength{\leftmargin}{\thmlistleft}%
			\setlength{\labelwidth}{\thmlistleft}%
			\setlength{\rightmargin}{\thmlistright}%
			\setlength{\partopsep}{\thmlistpartopsep}%
			\setlength{\topsep}{\thmlisttopsep}%
			\setlength{\parsep}{\thmlistparsep}%
			\setlength{\itemsep}{\thmlistitemsep}}}%
		    {\end{list}}%

\newcounter{prt}
\newenvironment{prt}{\begin{list}{\upshape (\alph{prt})}%
		    {\usecounter{prt}%
		        \setlength{\leftmargin}{\thmlistleft}%
			\setlength{\labelwidth}{\thmlistleft}%
			\setlength{\rightmargin}{\thmlistright}%
			\setlength{\partopsep}{\thmlistpartopsep}%
			\setlength{\topsep}{\thmlisttopsep}%
			\setlength{\parsep}{\thmlistparsep}%
			\setlength{\itemsep}{\thmlistitemsep}}}%
		    {\end{list}}%

\newcounter{rqm}
\newenvironment{rqm}{\begin{list}{\upshape (\arabic{rqm})}%
		    {\usecounter{rqm}%
		        \setlength{\leftmargin}{\thmlistleft}%
			\setlength{\labelwidth}{\thmlistleft}%
			\setlength{\rightmargin}{\thmlistright}%
			\setlength{\partopsep}{\thmlistpartopsep}%
			\setlength{\topsep}{\thmlisttopsep}%
			\setlength{\parsep}{\thmlistparsep}%
			\setlength{\itemsep}{\thmlistitemsep}}}%
		    {\end{list}}%

			{\end{list}}%

\newenvironment{itemlist}{\nopagebreak \begin{list}{$\bullet$}%
		       {\setlength{\leftmargin}{\thmlistleft}%
			\setlength{\labelwidth}{\thmlistleft}%
			\setlength{\rightmargin}{\thmlistright}%
			\setlength{\partopsep}{\thmlistpartopsep}%
			\setlength{\topsep}{\thmlisttopsep}%
			\setlength{\parsep}{\thmlistparsep}%
			\setlength{\itemsep}{\thmlistitemsep}}}%
			{\end{list}}%

			{\end{list}}%
\newlength{\myindent}
{\setlength{\myindent}{\parindent}\begin{list}{}%
			{\setlength{\leftmargin}{#1}\setlength{\rightmargin}{#1}%
			\setlength{\partopsep}{0pt}%
			\setlength{\topsep}{\thmtopspace}%
			\setlength{\parsep}{0pt}%
			\setlength{\itemsep}{0pt}}
			\item[]}
			{\end{list}}%

\newenvironment{proof*}{\begin{proof}}{ \end{proof}}

\newcommand{\dispand}[1][and]{\hbox to \hsize{#1 \hfill} \nonumber \\}

\newlength{\seqsplit}
\input{xypic}
\xyoption{all}

\title{Compactly generated homotopy categories}

\author{Henrik Holm \ \ \ and \ \ \ Peter J{\o}rgensen}

\address{\flushleft{Henrik} Holm, Department of Mathematical Sciences,
  University of Aarhus, Ny Munkegade, Building 530, DK--8000 Aarhus C,
  Denmark} 

\email{holm@imf.au.dk} 

\urladdr{http://home.imf.au.dk/holm}

\address{\flushleft{Peter} J{\o}rgensen, Department of Pure
  Mathematics, University of Leeds, Leeds LS2 9JT, United Kingdom}

\email{popjoerg@maths.leeds.ac.uk}

\urladdr{www.maths.leeds.ac.uk/\~{ }popjoerg}

\date{July 4, 2006}

\keywords{Compactly generated category, compact object,
  homotopy category, pure exact sequence, triangulated category}
 
\subjclass[2000]{16D20, 16D40, 16D50, 16D90, 16E05} 

\newcommand{\F}{\mathsf{F}}

\newcommand{\sK}{\mathsf{K}}
\newcommand{\sT}{\mathsf{T}}
\newcommand{\Hom}{\operatorname{Hom}}

\newcommand{\fpfunct}{\mathsf{f.p.funct}}

\newcommand{\X}{\mathbb{X}}

\newcommand{\Mod}{\mathsf{Mod}}
\newcommand{\fgmod}{\mathsf{mod}}
\newcommand{\Proj}{\mathsf{Proj}}
\newcommand{\Flat}{\mathsf{Flat}}
\newcommand{\Inj}{\mathsf{Inj}}

\newcommand{\PP}{\mathsf{PureProj}}
\newcommand{\PI}{\mathsf{PureInj}}

\begin{document}

\begin{abstract} 
  Over an associative ring we consider a class $\X$ of left modules
  which is closed under set-indexed coproducts and direct summands.
  We investigate when the triangulated homotopy category $\sK(\X)$ is
  compactly generated, and give a number of examples.
\end{abstract}

\maketitle
\setcounter{section}{-1}


\section*{Introduction}

\noindent
Let $\X$ be a class of left $R$--modules which is closed under
set-indexed coproducts and direct summands. When the (triangulated)
homotopy category $\sK(\X)$ associated to $\X$ is compactly generated,
it can be a powerful tool. Let us mention two recent examples from the
literature to illustrate this point:

\medskip
\noindent
In \cite[thm.\,2.4]{Jorgensen} it is proved that the homotopy category
$\sK(\Proj\,R)$ is compactly generated provided that $R$ is coherent
from either side, and that every flat left $R$--module has finite
projective dimension.

\medskip
\noindent
The result above is the cornerstone in proving that the class of
so-called Gorenstein projective modules is precovering (also called
contravariantly finite) when $R$ is commutative and noetherian with a
dualizing complex, cf.~\cite[cor.\,2.13]{Jorgensen2}. The question of
whether the Gorenstein projective modules really do constitute a
precovering class has been studied by many people; see for example
\cite{AvramovMartsinkovsky,EnochsRHA,Holm,Takahashi,EnochsJenda1}.

\medskip
\noindent
In \cite[prop.\,2.3]{Krause} it is shown that $\sK(\Inj\,R)$ is
compactly generated when $R$ is left noetherian. And in
\cite[cor.\,5.5]{IyengarKrause} this result is applied to give a new
and interesting characterization of Gorenstein rings in terms of
(totally) acyclic complexes of injective modules.

\medskip
\noindent
In this paper we study the general question of when $\sK(\X)$ is
compactly generated? More precisely, we give a number of sufficient
conditions on $R$ and $\X$ which ensure that $\sK(\X)$ is compactly
generated, and our results generalize those of
\cite{Jorgensen,Krause}.  At this point it is worth mentioning that
the innocent looking $\sK(\Mod\,\mathbb{Z})$ is known not to be
compactly generated, cf.~\cite[lem.\,E.3.2]{Neeman}.

\medskip
\noindent
Our main result is Theorem \eqref{thm:main} from Section
\ref{sec:main}. Sections \ref{sec:resolutions} and \ref{sec:Split}
develop the necessary machinery to provide us with examples where
Theorem \eqref{thm:main} can be applied. In the final Section
\ref{sec:summary} we use the previous results to list a number of
concrete classes $\X$ for which $\sK(\X)$ is compactly generated.


\section{Preliminaries} \label{sec:pure}

\noindent 
The assumptions, the notation, and the definitions from this section
will be used throughout the paper.

\begin{bfhpg}[Setup] \label{bfhpg:setup}
  Throughout, $R$ is a ring, and all modules are left $R$--modules
  unless otherwise specified. We use $R^\mathrm{op}$ to denote the
  opposite ring of $R$, and a left $R^\mathrm{op}$--module is
  naturally identified with a right $R$--module.  
  
  \medskip
  \noindent
  The symbol $\X$ always denotes a class of modules with
  \mbox{$\mathsf{Add}\,\X=\X$}, cf.~\eqref{bfhpg:notation} below.  For
  an arbitrary class of modules we write $\mathbb{A}$.
\end{bfhpg}

\begin{bfhpg}[Notation] \label{bfhpg:notation}
  We shall frequently use the following categories:
  \begin{itemlist}
  \item $\Mod\,R$ is the category of all $R$--modules, and $\fgmod\,R$
    is the category of all finitely presented $R$--modules.
  \item $\Proj\,R$, $\Inj\,R$, and $\Flat\,R$ are the categories of
    projective, injective, and flat $R$--modules, respectively.
  \item $\PP\,R$ and $\PI\,R$ are the categories of pure projective
    and pure injective $R$--modules, respectively,
    cf.~\cite[app.\,A]{JensenLenzing}.
  \item For a class of modules $\mathbb{A}$ we write
    $\mathsf{add}\,\mathbb{A}$ for the category of modules which are
    isomorphic to a direct summand of a module of the form $\coprod_{i
      \in I}A_i$, where \mbox{$A_i \in \mathbb{A}$}, and $I$ is a
    finite set. Allowing arbitrary index sets $I$ in this construction
    we get $\mathsf{Add}\,\mathbb{A}$.
  \end{itemlist} 
\end{bfhpg}

\begin{bfhpg}[The homotopy category]
  Let $\mathbb{A}$ be a class of modules with
  \mbox{$\mathbb{A}=\mathsf{add}\,\mathbb{A}$}. The objects of the
  \emph{homotopy category} $\sK(\mathbb{A})$ are chain complexes of
  modules from $\mathbb{A}$. Even though $\mathbb{A}$ is not abelian
  the notion of complexes is still well-defined since the condition
  \mbox{$\partial^2=0$} makes sense.  The morphisms of
  $\sK(\mathbb{A})$ are chain maps modulo homotopy equivalence.  

  \medskip
  \noindent
  By for example \cite[chap.\,10]{Weibel}, $\sK(\mathbb{A})$ carries
  the structure of a triangulated category with finite coproducts. If
  $\mathbb{A}$ has arbitrary set-indexed coproducts, then so has
  $\sK(\mathbb{A})$.
\end{bfhpg}

\begin{dfn} \label{dfn:compact}
  Let $\sT$ be a triangulated category, cf.~\cite{Neeman}, closed
  under set-indexed coproducts. An object \mbox{$C \in \sT$} is
  \emph{compact} if the natural map
  \begin{displaymath}
    \textstyle{\coprod}_{i \in I} \Hom_{\sT}(C,X_i) \longrightarrow
    \Hom_{\sT}(C,\textstyle{\coprod}_{i \in I} X_i) 
  \end{displaymath}
  is an isomorphism for any family $\{X_i\}_{i\in I}$ of objects in
  $\sT$. A set of objects \mbox{$\mathcal{G} \subseteq \sT$} is
  called \emph{a generating set} if the implication
  \begin{displaymath}
    \Hom_{\sT}(G,X) = 0 \ \text{ for all } \ G \in \mathcal{G} 
    \quad \Longrightarrow \quad X \cong 0
  \end{displaymath}
  holds for all \mbox{$X \in \sT$}. If $\sT$ has a generating set
  consisting of compact objects then $\sT$ is called \emph{compactly
    generated}.
\end{dfn} 

\begin{exa}
  The derived category $\mathsf{D}(\Mod\,R)$ of the abelian category
  $\Mod\,R$ is always compactly generated by the set
  \begin{displaymath}
    \mathcal{G} = \{\Sigma^nR\,|\,n\in \mathbb{Z}\}.
  \end{displaymath}
  Here $R$ is considered as a complex concentrated in degree zero with
  zero differentials, and $\Sigma^n \colon \mathsf{D}(\Mod\,R)
  \longrightarrow \mathsf{D}(\Mod\,R)$ denotes the $n$'th shift ``to
  the left'', that is, for a complex
  \begin{displaymath}
    X = \cdots \longrightarrow X_{s+1}
  \stackrel{\partial^X_{s+1}}{\longrightarrow} X_{s}
  \stackrel{\partial^X_{s}}{\longrightarrow} X_{s-1} \longrightarrow \cdots,
  \end{displaymath}
  the complex $\Sigma^nX$ has the module $X_{t-n}$ in degree $t$ and
  $(-1)^n\partial^X_{t-n}$ as its $t$'th differential.  

  \medskip
  \noindent
  Surprisingly, the corresponding homotopy category $\sK(\Mod\,R)$ is
  not even compactly generated when \mbox{$R=\mathbb{Z}$}; see
  \cite[lem.\,E.3.2]{Neeman}.  The ring $\mathbb{Z}$ has pure global
  dimension $1$, cf.~\eqref{bfhpg:pgldim}(1). It is a consequence of
  the theory developed in this paper, cf.~Section
  \ref{sec:summary}(3), that $\sK(\Mod\,R)$ is compactly generated
  when $R$ has pure global dimension $0$.
\end{exa}

\begin{rmk}
  For the notion of compact, that is, $\aleph_0$--compact objects the
  reader is referred to Neeman \cite[chap.\,4.1 and 4.2]{Neeman}. The
  definition of a genera\-ting set is taken from
  \cite[def.\,8.1.1]{Neeman}. Definition \eqref{dfn:compact} only has
  interest for us in the case where \mbox{$\sT = \sK(\X)$} for some
  class of modules $\X$, cf.~Setup \eqref{bfhpg:setup}.
\end{rmk}

\section{A condition for compact generatedness} \label{sec:main}

\noindent
In this section we give conditions on the module category $\X$,
cf.~Setup \eqref{bfhpg:setup}, which ensure that the associated
homotopy category $\sK(\X)$ is compactly generated. We begin by
stating our main result, but we postpone the proof until the end of
the section where the necessary preparations have been made.

\begin{thm} \label{thm:main}
  Let $\X=\mathsf{Add}\,\X$ be a class of $R$--modules, and assume
  that every finitely presented module $M$ has a right
  $\X$--resolution $X(M)$.  Then
  \begin{displaymath}
    \mathcal{G}_\X = \big\{ \Sigma^nX(M) \,\big|\,
    M\!\in\!\fgmod\,R\,, \
    n\in\mathbb{Z} \big\} 
  \end{displaymath}
  is a set of compact objects in $\sK(\X)$. Furthermore,
  $\mathcal{G}_\X$ generates $\sK(\X)$ if and only if every pure exact
  sequence in $\Mod\,R$, consisting of modules from $\X$, is split
  exact.
\end{thm}

\begin{bfhpg}[How to apply Theorem \eqref{thm:main}]
  In order to apply the theorem above we need examples of classes
  $\X=\mathsf{Add}\,\X$ satisfying:
  \begin{prt}
  \item Every finitely presented module has a right $\X$--resolution.
  \item Every pure exact sequence in $\Mod\,R$, consisting of modules
    from $\X$, is split exact.
  \end{prt}
  In Section \ref{sec:resolutions} we give examples of classes
  satisfying condition (a), and in Section \ref{sec:Split} we discuss
  how to check if (b) holds. In the final Section \ref{sec:summary} we
  use the results from the previous sections to list a number of
  concrete classes $\X$ for which $\sK(\X)$ is compactly generated.
\end{bfhpg}

\noindent
Before proving Theorem \eqref{thm:main} we need some preparation, in
particular some remarks about right resolutions and pure exactness.

\begin{bfhpg}[Right resolutions] \label{bfhpg:rightresolutions}
  Let $\mathbb{A}$ be a class of modules. An
  \emph{$\mathbb{A}$--pre\-en\-ve\-lo\-pe} of a module $M$ is a
  homomorphism $\varphi \colon M \longrightarrow A$ with $A \in
  \mathbb{A}$, such that given any other homomorphism $\varphi' \colon
  M \longrightarrow A'$ with $A' \in \mathbb{A}$ there is a (not
  necessarily unique) factorization,
  \begin{displaymath}
    \xymatrix{M \ar[r]^-{\varphi} \ar[d]_-{\varphi'} & A. \ar@{-->}[dl]
      \\ A' & {}
    }
  \end{displaymath}
  Let $M$ be a module which has an $\mathbb{A}$--preenvelope,
  \mbox{$\varphi^0 \colon M \longrightarrow A^0$}. Suppose that the
  cokernel of this map, $\mathrm{Coker}(\varphi^0)$, also has an
  $\mathbb{A}$--preenve\-lo\-pe, \mbox{$\varphi^1 \colon
    \mathrm{Coker}(\varphi^0) \longrightarrow A^1$}. If also
  $\mathrm{Coker}(\varphi^1)$ has an $\mathbb{A}$--preenvelope etc.,
  we may construct an \emph{augmented right $\mathbb{A}$--resolution}
  of $M$,
  \begin{displaymath}
    \tag{\text{$*$}}
    0 \longrightarrow M \stackrel{\varepsilon}{\longrightarrow} A^0
    \stackrel{\partial^0}{\longrightarrow} A^1
    \stackrel{\partial^1}{\longrightarrow} A^2
    \stackrel{\partial^2}{\longrightarrow} 
    \cdots.  
  \end{displaymath}
  Here $\varepsilon = \varphi^0$, and $\partial^n$ is the composition
  \begin{displaymath}
    \xymatrix{A^n \ar@{->>}[r] & \mathrm{Coker}(\varphi^n)
      \ar[r]^-{\varphi^{n+1}} & A^{n+1}.} 
  \end{displaymath}
  The complex $(*)$ is not necessarily exact (as
  $\mathbb{A}$--preenvelopes are not necessarily injective), however,
  $\Hom_R((*),A')$ is exact for every \mbox{$A' \in \mathbb{A}$}.
  There is a useful equivalent way of stating this property of the
  complex $(*)$, namely if we consider the chain map
  \begin{displaymath}
    \xymatrix@C=4ex{M \ar[d]^-{\varepsilon} & = & \cdots \ar[r] & 0
      \ar[d]^-{0} \ar[r] & 0 \ar[d]^-{0} \ar[r] & M 
      \ar[d]^-{\varepsilon} \ar[r] & 0 \ar[d]^-{0} \ar[r] & 0
      \ar[d]^-{0} \ar[r] & \cdots\\ 
    A(M) & = & \cdots \ar[r] & 0 \ar[r] & 0 \ar[r] & A^0
    \ar[r]^-{\partial^0} & A^1 \ar[r]^-{\partial^1} & A^2 \ar[r] &
    \cdots 
    }    
  \end{displaymath}
  then $\Hom_R(\varepsilon,A')$ is a quasi--isomorphism for all $A'
  \in \mathbb{A}$. In the given situation we refer to
  \begin{displaymath}
    A(M) \, = \ \ 0 \longrightarrow A^0
    \stackrel{\partial^0}{\longrightarrow} A^1
    \stackrel{\partial^1}{\longrightarrow} A^2
    \stackrel{\partial^2}{\longrightarrow} 
    \cdots
  \end{displaymath}
  as a (non-augmented) \emph{right $\mathbb{A}$--resolution} of $M$.
  Finally, the class $\mathbb{A}$ is called \emph{preenveloping} (also
  known as \emph{covariantly finite}) if every module has an
  $\mathbb{A}$--preenvelope, and thus a right
  $\mathbb{A}$--resolution.
  
  \medskip
  \noindent
  If a module $M$ admits a right $\mathbb{A}$--resolution it is in
  general not unique. However, by e.g.~\cite[ex.\,8.1.3]{EnochsRHA}
  all possible choices of right $\mathbb{A}$--resolutions $A(M)$ of
  $M$ are homotopy equivalent, and thus isomorphic in
  $\sK(\mathbb{A})$. Hence $A(M)$ is a well-defined object in the
  homotopy category of $\mathbb{A}$.  As the map \mbox{$\varepsilon
    \colon M \longrightarrow A(M)$} becomes a quasi--isomorphism
  whenever the functor $\Hom_R(-,A')$ is applied to it for $A' \in
  \mathbb{A}$, and since both $M$ and $A(M)$ are left-bounded
  complexes, \cite[prop.\,(2.7)(b)]{LWCAFHH} implies that
  \begin{displaymath}
    \xymatrix{\Hom_R(A(M),A')
      \ar[rrr]^-{\Hom_R(\varepsilon,A')}_-{\simeq} & & & \Hom_R(M,A')} 
  \end{displaymath}
  is a quasi--isomorphism for every \textsl{complex} $A'$ consisting
  of modules from $\mathbb{A}$. In particular we have an equivalence
  of functors $\sK(\mathbb{A}) \longrightarrow \Mod\,\mathbb{Z}$,
  \begin{displaymath}
    \tag{\text{$\dagger$}}
    \operatorname{H}_0\Hom_R(A(M),-) \,\simeq\,
    \operatorname{H}_0\Hom_R(M,-). 
  \end{displaymath}
  For the discussion above we do not assume that $\mathbb{A}$ has
  set-indexed coproducts, however, since we are interested in compact
  objects in the homotopy category we will focus on the case where
  $\mathbb{A}=\X$, cf.~\eqref{bfhpg:setup}.
\end{bfhpg}

\begin{prp} \label{prp:compact}
  If a module $M$ admits a right $\X$--resolution $X(M)$ then there is
  an equivalence of functors $\sK(\X) \longrightarrow \Mod\,\mathbb{Z}$,
  \begin{displaymath}
    \Hom_{\sK(\X)}(X(M),-) \,\simeq\,
    \operatorname{H}_0\Hom_R(M,-). 
  \end{displaymath}
  In particular, if $M$ is finitely generated and admits a right
  $\X$--resolution, then $X(M)$ is a compact object in $\sK(\X)$.
\end{prp}

\begin{proof}
  We have natural equivalences of functors,
  \begin{align*}
    \Hom_{\sK(\X)}(X(M),-) 
     &\,\simeq\,
    \operatorname{H}_0\Hom_R(X(M),-) \\ 
     &\,\simeq\,
    \operatorname{H}_0\Hom_R(M,-), 
  \end{align*}
  where the first isomorphism is standard, and the second is by
  $(\dagger)$ above. For the last claim in the proposition we use that
  $\operatorname{H}_0\Hom_R(M,-)$ commutes with set-indexed coproducts
  if $M$ is finitely generated.
\end{proof}

\begin{bfhpg}[Pure exact sequences] \label{bfhpg:pure}
  A (not necessarily short) sequence $Y$ in $\Mod\,R$ is \emph{pure
    exact} if and only if $\Hom_R(M,Y)$ is exact for all finitely
  presented modules $M$, cf.~\cite[thm.\,6.4]{JensenLenzing}.

\end{bfhpg}



\begin{proof}[Proof of Theorem \eqref{thm:main}]
  By Proposition \eqref{prp:compact} the set $\mathcal{G}_\X$ consists
  of compact objects. Strictly speaking $\mathcal{G}_\X$ is not a set,
  as $\fgmod\,R$ is not. However, we may of course restrict ourselves
  to just looking at isomorphism classes in $\fgmod\,R$, and they do
  constitute a set.

  \medskip
  \noindent
  Now, let $Y$ be an arbitrary object in $\sK(\X)$, that is, a chain
  complex of modules from $\X$. We claim that the following conditions
  are equivalent:
  \begin{eqc}
  \item $\Hom_{\sK(\X)}(\Sigma^nX(M),Y)=0$ for all $M \in \fgmod\,R$
    and $n \in \mathbb{Z}$;
  \item $Y$ is a pure exact sequence in $\Mod\,R$.
  \end{eqc}
  Having proved this, the last part of the theorem follows
  immediately, since an object \mbox{$Y \in \sK(\X)$} is isomorphic to
  zero if and only if $Y$ splits, cf.~\cite[ex.\,1.4.3]{Weibel}. The
  proof of the equivalence of $(i)$ and $(ii)$ follows from
  \eqref{bfhpg:pure} compared with the following calculation:
  \begin{align*}
    \Hom_{\sK(\X)}(\Sigma^nX(M),Y) &\cong
    \Hom_{\sK(\X)}(X(M),\Sigma^{-n}Y) \\
    &\cong \operatorname{H}_0\Hom_R(M,\Sigma^{-n}Y) \\
    &\cong \operatorname{H}_n\Hom_R(M,Y),
  \end{align*}  
  where the second isomorphism is by Proposition \eqref{prp:compact}.
\end{proof}

\section{Existence of right resolutions} \label{sec:resolutions}

\noindent
In this section we study constructions and examples of module classes
$\X$, cf.~Setup \eqref{bfhpg:setup}, for which every finitely
presented module has a right $\X$--resolution,
cf.~\eqref{bfhpg:rightresolutions}. This is of interest when we want
to apply our main Theorem \eqref{thm:main}. For reasons which will
become clear in Proposition \eqref{prp:PPandPI} and Example
\eqref{exa:pure} we will only focus on such classes which have the
additional property that they are contained in either $\PI\,R$ or
$\PP\,R$.

\medskip
\noindent
We begin by stating all our constructions and examples
\eqref{exa:inj}--\eqref{bfhpg:apply}, but we postpone the arguments to
the end of the section.

\begin{exa} \label{exa:inj}
  The following two examples are classical:
  \begin{prt}
  \item The class $\Inj\,R$ is preenveloping by \cite[thm.\,3.13 and
    3.26]{Rotman}.  However, in order for $\Inj\,R$ to be closed under
    coproducts, $R$ must be left noetherian,
    cf.~\cite[thm.\,4.27]{Rotman}.
  \item The class $\PI\,R$ is preenveloping by
    \cite[prop.\,7.6]{JensenLenzing}. However, in order for $\PI\,R$
    to be closed under coproducts, $R$ must be pure-semisimple, and in
    this case we actually have $\PI\,R = \Mod\,R$,
    cf.~\cite[thm.\,B.18]{JensenLenzing}.
  \end{prt} 
\end{exa}

\noindent
Actually, Example \eqref{exa:inj}(a) admits the following
generalization.  In the result below, $\mathsf{Inj}(M_1,\ldots,M_n)$
is defined as in \eqref{dfn:DefOfMinj}.

\begin{prp} \label{prp:Minj}
  Assume that $R$ is left noetherian and let $M_1,\ldots,M_n$ be
  $R$--bimodules such that each $M_j$ is finitely generated as a left
  $R$--module.  Then every module has a right
  $\mathsf{Inj}(M_1,\ldots,M_n)$--resolution.
\end{prp}

\noindent
Next we present examples which are contained in $\PP\,R$:

\begin{exa} \label{exa:PROJ}
  The following conclusions hold:
  \begin{prt}
  \item Every finitely presented module has a right resolution with
    respect to $\PP\,R$.
  \item If $R$ is right coherent then every finitely presented module
    has a right resolution with respect to
    \mbox{$\mathsf{Add}\,({}_RR) = \Proj\,R$}.
  \end{prt}
\end{exa} 

\noindent
Actually, Example \eqref{exa:PROJ}(b) admits a generalization:

\begin{prp} \label{prp:AddM}
  Assume that $R$ is right coherent, and let $M$ be an $R$--bimodule
  which is finitely presented from either side.  Then every finitely
  presented module has a right $\mathsf{Add}({}_RM)$--resolution.
\end{prp}

\noindent
In the next result, \mbox{$\mathsf{Gproj}\,R$} denotes the class of
finitely generated Gorenstein projective modules,
cf.~\cite{Auslander,AuslanderBridger}, and \mbox{$\mathsf{GFlat}\,R$}
is the class of Gorenstein flat modules,
cf.~\cite{EnochsJendaTorrecillas}. Furthermore,
$\varinjlim\mathsf{Gproj}\,R$ is the class of modules which can be
written as a colimit in $\Mod\,R$ of some functor \mbox{$I
  \longrightarrow \mathsf{Gproj}\,R$}, where $I$ is a small filtering
category, cf.~\eqref{bfhpg:support}.

\begin{prp} \label{prp:Gproj}
  Assume that $R$ is commutative and noetherian with a dualizing
  complex. If \mbox{$\,\varinjlim\mathsf{Gproj}\,R=\mathsf{GFlat}\,R$}
  (this happens for example if $R$, in addition, has finite injective
  dimension over itself) then every fi\-ni\-te\-ly presented module
  has a right $\mathsf{Add}(\mathsf{Gproj}\,R)$--resolution
\end{prp}

\noindent
In the next result, $\mathsf{sub}\,\mathbb{A}$ is defined as in
Lemma \eqref{lem:sub}.

\begin{prp} \label{prp:Addsub}
  Let $\mathbb{A}$ be a class of modules, and assume that every
  finitely presented module $N$ has an $\mathbb{A}$--preenvelope
  \mbox{$\varphi \colon N \longrightarrow A$} for which the image
  $\operatorname{Im}(\varphi)$ is finitely presented. Then every
  finitely presented module has a right
  $\mathsf{Add}(\mathsf{sub}\,\mathbb{A})$--resolution.
\end{prp}

\begin{bfhpg}[How to apply Proposition \eqref{prp:Addsub}] \label{bfhpg:apply}
  To apply the proposition above we first of all need a class
  $\mathbb{A}$ such that every finitely presented module has an
  $\mathbb{A}$--preenvelope. For instance, $\mathbb{A}$ could be any
  of the following preenveloping classes:
  \begin{prt}
  \item \mbox{$\PI\,R \cap \Flat\,R$} if $R$ is right coherent,
    cf.~\cite[prop.\,6.6.6]{EnochsRHA}.
  \item $\Flat\,R$ if $R$ is right coherent,
    cf.~\cite[prop.\,6.5.1]{EnochsRHA}.
  \item The class of $S$--torsion free modules, when $R$ is
    commutative and \mbox{$S \subseteq R$} is a multipli\-cative
    subset, cf.~\eqref{bfhpg:tor} below.
  \end{prt}
  But other choices of $\mathbb{A}$ are also possible; for example
  from the proof of Proposition \eqref{prp:AddM} it will follow that:
  \begin{prt}
  \setcounter{prt}{3}
\item If $R$ is right coherent, and $M$ is an $R$--bimodule which is
  finitely presented from either side, then every finitely presented
  module has an $\mathsf{add}({}_RM)$--preenvelope (note
  $\mathsf{add}$, not $\mathsf{Add}$).
  \end{prt}
  
  \noindent
  However, $\mathbb{A}$ must have the additional property that among
  all the preenvelopes of a given finitely presented module $N$, there
  should exist one with a finitely presented image. We note that
  \begin{rqm}
  \item If $R$ is left noetherian then the image of every
    $\mathbb{A}$--preenvelope of $N$ is finitely presented.
  \item If \mbox{$\mathbb{A} \subseteq \fgmod\,R$} (this is the case
    in (d) above) and $R$ is left coherent, then the image of every
    pre\-envelope of $N$ is finitely presented by
    \cite[thm.\,3.2.24]{EnochsRHA}.
  \end{rqm}
\end{bfhpg}

\noindent
Before proving \eqref{prp:Minj}---\eqref{prp:Addsub} we will get
\eqref{bfhpg:apply}(c) out of the way:

\begin{bfhpg}[$S$--torsion free modules] \label{bfhpg:tor}
  Let $R$ be commutative and let \mbox{$S \subseteq R$} be a
  multiplicative subset. For any module $M$ its \emph{$S$--torsion
    submodule} is defined as
  \begin{displaymath}
    \Gamma_SM = \big\{ x \in M \,\big|\, sx=0 \text{ for some }
    s \in S \big\}.
  \end{displaymath}
  We say that $M$ is \emph{$S$--torsion free} if \mbox{$\Gamma_SM=0$}.
  It is easy to see that the class of $S$--torsion free modules
  preenveloping as $M \longrightarrow M/\Gamma_SM$ is an $S$--torsion
  free preenevelope of $M$.
\end{bfhpg}



\noindent
In the rest of the section we prove
\eqref{prp:Minj}---\eqref{prp:Addsub}. We begin with a

\begin{dfn} \label{dfn:DefOfMinj}
  Assume that $R$ is left noetherian and let $M_1,\ldots,M_n$ be
  $R$--bimodules such that each $M_j$ is finitely generated as a left
  $R$--module.  A module $J$ belongs to $\mathsf{Inj}(M_1,\ldots,M_n)$
  if and only if there exist injective modules $I_1,\ldots,I_n$ such
  that $J$ is a direct summand of
  \begin{displaymath}
    \Hom_R(M_1,I_1) \oplus \cdots \oplus \Hom_R(M_n,I_n).
  \end{displaymath}
\end{dfn}

\begin{proof}[Proof of Proposition \eqref{prp:Minj}]
  Since $R$ is left noetherian there is by \cite[proof of
  thm.\,5.4.1]{EnochsRHA} a set of injective modules $\mathbb{E}$ such
  that \mbox{$\Inj\,R = \mathsf{Add}\,\mathbb{E}$}. Using that
  $\Hom_R(M_j,-)$ commutes with set-indexed coproducts we see that
  $\mathsf{Inj}(M_1,\ldots,M_n)$ has the form
  $\mathsf{Add}\,\mathbb{A}$, where
  \begin{displaymath}
    \mathbb{A} = \big\{\Hom_R(M_j,E)
    \, \big|\ 
    j\in\{1,\ldots,n\} \,,\, E \in \mathbb{E} \big\}.
  \end{displaymath}
  From the description in \eqref{dfn:DefOfMinj} it is clear that
  $\mathsf{Inj}(M_1,\ldots,M_n)$ is closed under products, and hence
  it follows easily from \cite[prop.\,6.2.1]{EnochsRHA} that
  $\mathsf{Inj}(M_1,\ldots,M_n)$ is preenveloping.
\end{proof}

\begin{rmk}
  $\mathsf{Inj}(M_1,\ldots,M_n)$ is contained in $\PI\,R$.
\end{rmk}

\begin{proof}[Proof of Example \eqref{exa:PROJ}]
  Part (a) is clear as $\PP\,R$ contains every finitely presented
  module. Part (b) follows from \cite[exa.\,3.4]{EJbalanced}.
\end{proof}

\noindent
Before we go on we need a few facts about finitely presented modules:

\begin{lem} \label{lem:quotient}
  The following conclusions hold:
  \begin{prt}
  \item If $M$ is a finitely presented module and \mbox{$S \subseteq\!
      M$} a finitely gene\-ra\-ted submodule then the quotient $M/S$
    is finitely presented.
  \item If $M$ is finitely generated and \mbox{$S \subseteq\! M$} is a
    submodule such that $M/S$ is finitely presented, then $S$ is
    finitely generated.
  \item Assume that $R$ is left coherent, $M$ is finitely generated,
    and $N$ is finitely presented. If \mbox{$\varphi \colon M
      \longrightarrow N$} is a homomorphism then
    $\mathrm{Ker}(\varphi)$ is finitely generated.
  \end{prt}
\end{lem}

\begin{proof}
  Part (a) is easy to prove, and part (b) can be found in for
  ex\-am\-ple \cite[prop.\,(4.26)(b)]{Lam}. Using (b) we can easily
  prove (c):

  \medskip
  \noindent
  As $M$ is finitely generated then so is $\mathrm{Im}(\varphi)$.  As
  $R$ is coherent and $\mathrm{Im}(\varphi)$ is a finitely generated
  submodule of the finitely presented module $N$, it follows by
  \cite[def.\,(4.51) and cor.\,(4.52)]{Lam} that
  $\mathrm{Im}(\varphi)$ is even finitely presented. Applying (b) to
  the inclusion \mbox{$\mathrm{Ker}(\varphi) \subseteq\! M$}, which
  has \mbox{$M/\mathrm{Ker}(\varphi)\cong \mathrm{Im}(\varphi)$}, we
  get that $\mathrm{Ker}(\varphi)$ is finitely generated. \qedhere
\end{proof}

\begin{lem} \label{lem:Add-resolution}
  Assume that \mbox{$\mathbb{A} \subseteq \fgmod\,R$}, and that every
  finitely presented module has an $\mathbb{A}$--preenvelope. Then
  every finitely presented module has a right
  $\mathsf{Add}\,\mathbb{A}$--resolution.
\end{lem} 

\begin{proof}
  Let $M$ be a finitely presented module, and let \mbox{$\varphi^0
    \colon M \longrightarrow A^0$} be an $\mathbb{A}$--preenvelope.
  Since $\mathrm{Im}(\varphi^0)$ is finitely generated and $A^0$ is
  finitely presented, Lemma \eqref{lem:quotient}(a) implies that
  $\mathrm{Coker}(\varphi^0)$ is finitely presented, so it has an
  $\mathbb{A}$--preenvelope, \mbox{$\varphi^1 \colon
    \mathrm{Coker}(\varphi^0) \longrightarrow A^1$}. Conti\-nuing in this
  manner we build an augmented right $\mathbb{A}$--resolution of $M$,
  \begin{displaymath}
    A^+(M) = \, 0 \longrightarrow M \longrightarrow A^0
    \longrightarrow A^1 \longrightarrow \cdots.
  \end{displaymath}
  To finish the proof it suffices to see that $\Hom_R(A^+(M),A')$ is
  exact for every \mbox{$A' \in \mathsf{Add}\,\mathbb{A}$}. We may
  assume that $A'$ has the form \mbox{$\coprod_{i \in I}A_i$} where
  \mbox{$A_i \in \mathbb{A}$}. Finally we simply have to use that
  \begin{displaymath}
    \Hom_R(A^+(M),\textstyle{\coprod}_{i \in I}A_i) \,\cong\,
    \textstyle{\coprod}_{i \in I}\Hom_R(A^+(M),A_i) 
  \end{displaymath}
  as every module in $A^+(M)$ is, in particular, finitely generated.
\end{proof}

\begin{bfhpg}[Modules with support in a category] \label{bfhpg:support}
  For a class \mbox{$\mathbb{A} = \mathsf{add}\,\mathbb{A}$} of
  finitely presented modules, Lenzing \cite{Lenzing} studies the class
  \mbox{$\varinjlim \mathbb{A}$} of all colimits in $\Mod\,R$ of
  functors \mbox{$I \longrightarrow \mathbb{A}$}, where $I$ is a small
  filtering category, cf.~\cite[chap.\,IX]{MacLane}.

  \medskip
  \noindent
  A module in \mbox{$\varinjlim \mathbb{A}$} is said to have
  \emph{support} in $\mathbb{A}$, and \cite[prop.\,2.1]{Lenzing} gives
  two other characterizations of these modules.  The following result
  can be found in for example \cite[thm.\,3.2]{Angeleri-Hugel} or
  \cite[sec.\,(4.2)]{C-B}: 
\end{bfhpg}

\begin{prp} \label{prp:flat-preenveloping}
  If \mbox{$\mathbb{A} \subseteq \fgmod\,R$} with
  \mbox{$\mathsf{add}\,\mathbb{A} = \mathbb{A}$}, then the following
  two conditions are equivalent:
  \begin{eqc}
  \item Every finitely presented module has an
    $\mathbb{A}$--preenvelope.
  \item \mbox{$\varinjlim \mathbb{A}$} is closed under products.
  \end{eqc}
\end{prp}

\noindent
The following proof is a consequence of the proposition above:

\begin{proof}[Proof of Proposition \eqref{prp:Gproj}]
  Since $R$ is commutative and noetherian with a dualizing complex,
  \cite[thm.\,(5.7)]{LWCAFHH} gives that $\mathsf{GFlat}\,R$ is closed
  under products. The assumption \mbox{$\varinjlim(\mathsf{Gproj}\,R)
    = \mathsf{GFlat}\,R$}, Proposition \eqref{prp:flat-preenveloping}
  and Lemma \eqref{lem:Add-resolution} applied to
  \mbox{$\mathbb{A}=\mathsf{Gproj}\,R$} give that every finitely
  presented module has a right resolution with respect to
  $\mathsf{Add}(\mathsf{Gproj}\,R)$. 

  \medskip
  \noindent
  It remains to prove the claim in parentheses, namely that the
  equality \mbox{$\varinjlim(\mathsf{Gproj}\,R) = \mathsf{GFlat}\,R$}
  holds when $R$ is commutative and noetherian with finite
  in\-jec\-tive dimension over itself. The inclusion ``$\supseteq$''
  follows from \cite[thm.\,10.3.8]{EnochsRHA}. The opposite inclusion
  ``$\subseteq$'' follows from combining \cite[thm.\,(3.5)]{LWCAFHH}
  and \cite[cor.\,2.6.17]{Weibel} with the fact that every module has
  finite Gorenstein flat dimension. For the latter claim see for
  example \cite[(1.3) and thm.\,(4.1)]{LWCAFHH}, or
  \cite[thm.\,(5.2.10)]{LWC} in the local case.
\end{proof}

\noindent
Our next goal is to provide the proof of Proposition \eqref{prp:AddM}:

\begin{proof}[Proof of Proposition \eqref{prp:AddM}]
  As ${}_RM$ is finitely presented it suffices by Lemma
  \eqref{lem:Add-resolution} to show that every finitely presented
  module $N$ has an \mbox{$\mathsf{add}({}_RM)$}--preenvelope.  We
  start by proving that the $R^\mathrm{op}$--module $\Hom_R(N,M)$ is
  finitely generated: Since $N$ is finitely presented there is an
  exact sequence,
  \begin{displaymath}
    F_1 \longrightarrow F_0 \longrightarrow N \longrightarrow 0,
  \end{displaymath}
  where \mbox{$F_i \cong ({}_RR)^{b_i}$} is finitely generated and
  free.  Applying the left exact functor $\Hom_R(-,M)$ to this
  sequence we get
  \begin{displaymath}
    \tag{\text{$*$}}
    0 \longrightarrow \Hom_R(N,M) \longrightarrow \Hom_R(F_0,M)
    \longrightarrow \Hom_R(F_1,M).
  \end{displaymath}
  Since \mbox{$\Hom_R(F_i,M) \cong (M_R)^{b_i}$}, and since $M_R$ is
  finitely presented, we see that $\Hom_R(F_i,M)$ is finitely
  presented.  Applying Lemma \eqref{lem:quotient}(c) to $(*)$ we get
  that $\Hom_R(N,M)$ is finitely generated, and we write
  \begin{displaymath}
    \tag{\text{$**$}}
    \Hom_R(N,M) = h_1R + \cdots + h_tR. 
  \end{displaymath}
  We claim that the map $\varphi \colon N \longrightarrow M^t$ defined
  by
  \begin{displaymath}
    z \,\longmapsto\, (h_1(z),\ldots,h_t(z)) 
  \end{displaymath}
  is an $\mathsf{add}({}_RM)$--preenvelope of $N$.  To see this it
  suffices to prove that any homomorphism \mbox{$\psi \colon N
    \longrightarrow M^k$} from $N$ to a finite power of ${}_RM$ lifts
  to $M^t$,
  \begin{displaymath}
    \xymatrix@R=3ex@C=3ex{
      N \ar[dr]_-{\psi} \ar[rr]^-{\varphi} & & M^t.
      \ar@{-->}[dl]^-{u} \\ 
      & M^k &  
    }
  \end{displaymath}
  Furthermore, without loss of generality we may assume that $k=1$.
  To define $u$ use $(**)$ to write $\psi \in \Hom_R(N,M)$ as
  \begin{displaymath}
    \psi =  h_1 r_1 \,+\, \cdots
    \,+\, h_t r_t
  \end{displaymath}
  for suitable $r_1,\ldots,r_t \in R$. We can then define $u \colon
  M^t \to M$ by
  \begin{displaymath}
    ( x_1,\ldots,x_t )  \,\longmapsto\,
    x_1 r_1 + \cdots + x_t r_t.
  \end{displaymath}
  Now \mbox{$u\varphi = \psi$} because for $z \in N$ we have:
  \begin{align*}
    u\varphi(z) &= u (h_1(z),\ldots,h_t(z)) \\
    &= h_1(z) r_1 \,+\, \cdots \,+\, h_t(z) r_t \\ 
    &= \psi(z). \qedhere
  \end{align*}
\end{proof}

\noindent
Finally we need to show Proposition \eqref{prp:Addsub}, but first a
little preparation:

\begin{lem} \label{lem:sub}
  Let $\mathbb{A}$ be any class of modules and define
  \begin{displaymath}
    \mathsf{sub}\,\mathbb{A} = \{ S\in \fgmod\,R \,|\,
    \text{$S \subseteq A$ for some $A \in \mathbb{A}$} \}.
  \end{displaymath}
  Assume that $M$ is a finitely presented module, and that $M$ admits
  an $\mathbb{A}$--preenvelope $\varphi \colon M \longrightarrow A$
  such that $\operatorname{Im}(\varphi)$ is finitely presented. Then
  $M$ has a right $\mathsf{sub}\,\mathbb{A}$--preenvelope.
\end{lem}

\begin{proof}
  By assumption $M$ has an $\mathbb{A}$--preenvelope \mbox{$\varphi
    \colon M \longrightarrow A$} such that $\mathrm{Im}(\varphi)$ is
  finitely presented. Consider the obvious factorization,
  \begin{displaymath}
    \xymatrix@R=3ex@C=3ex{ M \ar@{->>}[dr]_-{\tilde{\varphi}}
      \ar[rr]^-{\varphi} & & A \\ & \mathrm{Im}(\varphi)
      \ar@{^(->}[ur]_-{i} & 
    }
  \end{displaymath}
  By definition the module $\mathrm{Im}(\varphi)$ belongs to
  $\mathsf{sub}\,\mathbb{A}$, and it is easy to verify that
  $\tilde{\varphi} \colon M \longrightarrow \mathrm{Im}(\varphi)$ is
  indeed a $\mathsf{sub}\,\mathbb{A}$--preenvelope of $M$.
\end{proof}

\begin{proof}[Proof of Proposition \eqref{prp:Addsub}]
  The assumptions on $\mathbb{A}$ and Lemma \eqref{lem:sub} ensure
  that every finitely presented module has a
  $\mathsf{sub}\,\mathbb{A}$--preenvelope.  Since
  $\mathsf{sub}\,\mathbb{A}$ is contained in $\fgmod\,R$ by
  definition, Lemma \eqref{lem:Add-resolution} finishes the proof.
\end{proof}

\begin{rmk}
  The class $\mathsf{Add}(\mathsf{sub}\,\mathbb{A})$ is contained in
  $\PP\,R$ since $\mathsf{sub}\,\mathbb{A}$ is contained in
  $\fgmod\,R$ by definition. 
\end{rmk}

\section{When does a pure exact sequence split?} \label{sec:Split}

\noindent
Given a class of modules $\mathbb{A}$, we discuss in this section how
to see if every pure exact sequence in $\Mod\,R$, consisting of
modules from $\mathbb{A}$, is split exact. This question is of
interest when we wish to apply our main Theorem \eqref{thm:main}. We
begin by outlining the idea of this section, but we postpone the
arguments until later:

\medskip
\noindent
For any class of modules we consider two conditions \textsc{(pp)} and
\textsc{(pi)}, cf.~Definition \eqref{dfn:pure}. These conditions can
be checked, and are indeed fulfilled in many cases as Example
\eqref{exa:pure} below shows. The conditions \textsc{(pp)} and
\textsc{(pi)} are the key ingredients in the following proposition,
which is the main result of this section:

\begin{prp} \label{prp:PPandPI}
  Let $\mathbb{A}$ be a class of modules satisfying at least one of
  the two conditions \textsc{(pp)} or \textsc{(pi)}.  Then every pure
  exact sequence in $\Mod\,R$, consisting of modules from
  $\mathbb{A}$, is split exact.
\end{prp}

\begin{exa} \label{exa:pure}
  The conclusions below hold:
  \begin{prt}
  \item If $R$ has finite left pure global dimension,
    cf.~\eqref{bfhpg:pgldim}, then every subclass of $\PP\,R$
    satisfies \textsc{(pp)}, and every subclass of $\PI\,R$ satisfies
    \textsc{(pi)}.
  \item If every flat $R$--module has finite projective dimension then
    every subclass of $\Proj\,R$ satisfies \textsc{(pp)}.
  \item If $R$ is left noetherian then every subclass of $\Inj\,R$
    satisfies \textsc{(pi)}.
  \end{prt}
\end{exa}

\noindent
The rest of the section is devoted to proving Proposition
\eqref{prp:PPandPI} and Example \eqref{exa:pure}. We begin with the
following:

\begin{dfn} \label{dfn:F}
  For a class of modules $\mathbb{A}$ we define $\F(\mathbb{A})$ to be
  the class of modules which are isomorphic to some kernel
  (equivalently, some image, or some cokernel) in a pure exact
  sequence,
  \begin{displaymath}
    \cdots \longrightarrow A_{n+1} \longrightarrow A_{n}
    \longrightarrow A_{n-1} \longrightarrow \cdots
  \end{displaymath}
  where every $A_n$ belongs to $\mathsf{Add}\,\mathbb{A}$.
\end{dfn}

\noindent
The properties \textsc{(pp)} and \textsc{(pi)} for a class
$\mathbb{A}$, which occur in this section's main result
\eqref{prp:PPandPI}, are defined in terms of $\F(\mathbb{A})$ from
Definition \eqref{dfn:F}:

\begin{dfn} \label{dfn:pure}
  For a class of modules $\mathbb{A}$ we consider the properties:
  \begin{prt}
  \item[(\textsc{pp})] There exists a \mbox{$d \geqslant 0$} such that
    for every \mbox{$M \in \F(\mathbb{A})$} and every pure exact
    sequence, $0 \to K_d \to A_{d-1} \to \cdots \to A_0 \to M \to 0$,
    with $A_0,\ldots,A_{d-1} \in \mathbb{A}$, the module $K_d$ must be
    pure projective.
  \item[(\textsc{pi})] There exists a \mbox{$d \geqslant 0$} such that
    for every \mbox{$M \in \F(\mathbb{A})$} and every pure exact
    sequence, $0 \to M \to A^0 \to \cdots \to A^{d-1} \to C^d \to 0$,
    with $A^0,\ldots,A^{d-1} \in \mathbb{A}$, the module $C^d$ must be
    pure injective.
  \end{prt}
\end{dfn}

\noindent
The purpose of the following Observation \eqref{obs:contained} and the
subsequent Lemma \eqref{lem:FProjandFInj} is to get a better feeling
for the construction $\F(-)$ from Definition \eqref{dfn:F}.

\begin{obs} \label{obs:contained}
  Clearly \mbox{$\mathsf{Add}\,\mathbb{A} \subseteq \F(\mathbb{A})$}.
  Furthermore, if \mbox{$\mathbb{A}=\mathsf{add}\,\mathbb{A}$}
  consists of finitely presented modules then \mbox{$\F(\mathbb{A})
    \subseteq \varinjlim \mathbb{A}$}, cf.~\eqref{bfhpg:support}, as:

  \medskip
  \noindent
  If \mbox{$M \in \F(\mathbb{A})$} then in particular there exist a
  module \mbox{$A \in \mathsf{Add}\,\mathbb{A}$} and a pure
  monomorphism \mbox{$0 \to M \to A$}. Since \mbox{$A \in \varinjlim
    \mathbb{A}$}, and since \mbox{$\varinjlim \mathbb{A}$} is closed
  under pure submodules by \cite[prop.\,2.2]{Lenzing}, it follows that
  $M$ belongs to \mbox{$\varinjlim \mathbb{A}$}.
\end{obs}

\begin{lem} \label{lem:FProjandFInj}
  The following conclusions hold:
  \begin{prt}
  \item \mbox{$\F(\Proj\,R) \subseteq \Flat\,R$}. If $R$ is a
    commutative integral domain which is not a field, then the
    inclusion is strict.
  \item If $R$ is left noetherian then $\F(\Inj\,R)=\Inj\,R$.
  \end{prt}
\end{lem}

\begin{proof}
  ``(a)'': It follows immediately from Definition \eqref{dfn:F} of
  $\F(-)$ that $\F(\Proj\,R) = \F(\mathsf{add}\,R)$, and the latter is
  contained in \mbox{$\varinjlim(\mathsf{add}\,R) = \Flat\,R$} by
  Observation \eqref{obs:contained}.  If $R$ is a commutative integral
  domain with quotient field \mbox{$Q \neq R$} then $Q$ belongs to
  \mbox{$\Flat\,R$}, but \mbox{$Q \notin \F(\Proj\,R)$} since $Q$
  cannot even be embedded into a free module.
  
  \medskip
  \noindent
  ``(b)'': Only the inclusion \mbox{$\F(\Inj\,R) \subseteq \Inj\,R$}
  is non-trivial; thus we let \mbox{$M \in \F(\Inj\,R)$} and use
  Baer's criterion to show that $M$ is injective: The assumption on
  $M$ implies, in particular, the existence of a pure epimorphism $f
  \colon I \twoheadrightarrow M$, where $I$ is injective.  Let
  \mbox{$\mathfrak{a} \subseteq R$} be an ideal, and let \mbox{$i
    \colon \mathfrak{a} \to R$} be the inclusion. Given a homomorphism
  \mbox{$u \colon \mathfrak{a} \to M$} we must find \mbox{$v \colon R
    \to M$} with \mbox{$vi = u$}.  Since $R$ is left noetherian the
  ideal $\mathfrak{a}$ is finitely presented, so by assumption on $f$
  we get $g \colon \mathfrak{a} \to I$ with $f g = u$,
  \begin{displaymath}
    \xymatrix@R=7ex@C=7ex{\mathfrak{a} \ar@{..>}[dr]^-{g} \ar[d]_-{u}
      \ar@{^(->}[r]^-{i} & R 
      \ar@{..>}[d]^-{h} \\ M & I \ar@{->>}[l]^-{f} }
  \end{displaymath}
  Injectivity of $I$ then gives \mbox{$h \colon R \to I$} with
  \mbox{$hi = g$}. Consequently, the map \mbox{$v=f h \colon
    R \to M$} is the desired one.
\end{proof}

\begin{bfhpg}[Krull dimension of categories] \label{bfhpg:Krull}
  Geigle \cite[def.\,2.1]{Geigle} has introduced a Krull dimension for
  a small additive category $\mathcal{C}$. By definition, the
  Krull--Geigle dimension of $\mathcal{C}$ coincides with the
  Krull--Gabriel dimension (introduced in \cite{Gabriel} using
  filtrations of localizing subcategories) of
  $\fpfunct(\mathcal{C}^\mathrm{op},\mathsf{Ab})$. The
  latter is the category of all covariant, additive, and finitely
  presented functors $\mathcal{C}^\mathrm{op} \longrightarrow
  \mathsf{Ab}$, where $\mathsf{Ab}=\mathsf{Mod}\,\mathbb{Z}$.

  \medskip
  \noindent
  For a ring $R$, Jensen--Lenzing \cite[pp.\,197--199]{JensenLenzing}
  consider a Krull dimension for $\mathsf{mod}(R^\mathrm{op})$; by
  definition it is the Krull--Gabriel dimension of the category
  $\fpfunct(\mathsf{mod}(R^\mathrm{op}),\mathsf{Ab})$.

  \medskip
  \noindent
  For an Artin algebra $\Lambda$ there is by
  \cite[thm.\,3.3]{AuslanderReitenSmalo} a duality, in other words a
  ``contravariant equivalence'',
  \begin{displaymath}
    D \colon \mathsf{mod}\,\Lambda \longrightarrow
    \mathsf{mod}(\Lambda^\mathrm{op}). 
  \end{displaymath}
  Consequently, there is also an equivalence of categories,
  \begin{displaymath}
    \xymatrix@R=0ex@C=3ex{
    \fpfunct(\mathsf{mod}(\Lambda^\mathrm{op}),\mathsf{Ab}) 
    \ar[rr]^-{\sim} & &  
    \fpfunct((\mathsf{mod}\,\Lambda)^\mathrm{op},\mathsf{Ab}) \\ 
    F \qquad \ar@{|->}[rr] & & \qquad F \circ D.
    }
  \end{displaymath}
  In particular, the Krull--Jensen--Lenzing dimension of
  $\mathsf{mod}(\Lambda^\mathrm{op})$ agrees with the Krull--Geigle
  dimension of $\mathsf{mod}\,\Lambda$. By \cite[thm.\,4.3]{Geigle}
  the latter is finite when $\Lambda$ is a tame hereditary Artin
  algebra.
\end{bfhpg}

\begin{bfhpg}[Pure global dimension] \label{bfhpg:pgldim}
  The definition of the \emph{left pure global dimension} for a ring
  $R$, denoted \mbox{$\operatorname{l.\!p.\!gl.\!dim}R$}, may be found
  in for example \cite[def.\,A.14]{JensenLenzing}.  Below we give
  examples of classes of rings with finite left pure global dimension:
  \begin{rqm}
  \item If $\max\{\aleph_0,|R|\}=\aleph_t$ then
    \mbox{$\operatorname{l.\!p.\!gl.\!dim}R\leqslant t+1$} by
    \cite[\S2]{GrusonJensen}.

  \item If \mbox{$I \subseteq R$} is a two-sided ideal then
    \mbox{$\operatorname{l.\!p.\!gl.\!dim}R/I\leqslant
      \operatorname{l.\!p.\!gl.\!dim}R$}. If $R$ is commutative and
    \mbox{$S \subseteq R$} is a multiplicative subset then
    \mbox{$\operatorname{l.\!p.\!gl.\!dim}S^{-1}R\leqslant
      \operatorname{l.\!p.\!gl.\!dim}R$},
    cf.~\cite[prop.\,1.1]{KielpinskiSimon}.
    
  \item If $(R,\mathfrak{m},k)$ is a commutative local noetherian
    domain of Krull dimension $1$, and $k$ is at most countable, then
    \mbox{$\operatorname{l.\!p.\!gl.\!dim}R=1$} by
    \cite[prop.\,4.7]{JensenLenzing2}.
    
  \item If $R$ is a finite dimensional $k$--algebra ($k$ any field) of
    tame representation type, which is either hereditary or a
    radical-squared zero algebra, then
    \mbox{$\operatorname{l.\!p.\!gl.\!dim}R\leqslant 2$} by
    \cite[prop.\,3.3]{Baer}.
    
  \item Some specific examples of four dimensional $k$--algebras
    (which are neither hereditary, nor radical-squared zero) with
    finite pure global dimension may be found in
    \cite[prop.\,5.1]{Baer}. The reader might also want to consult
    \cite[cor.\,11.33 and 11.34]{JensenLenzing}.
  
  \item If the category $\mathsf{mod}(R^\mathrm{op})$ has finite Krull
    dimension $d$ according to Jensen--Lenzing
    \cite[pp.\,197--199]{JensenLenzing}, cf.~\eqref{bfhpg:Krull}
    above, then \mbox{$\operatorname{l.\!p.\!gl.\!dim}R\leqslant d$}
    by \cite[thm.\,11.31]{JensenLenzing}. This applies for example to
    the rings:
    \begin{itemlist}
    \item $R$ is a Dedekind domain, cf.~\cite[thm.\,8.55,
      cor.\,11.32]{JensenLenzing}.
    \item $R$ is a tame hereditary Artin algebra,
      cf.~\eqref{bfhpg:Krull}.
    \end{itemlist}
    
  \item If $R$ if von Neumann regular then every exact sequence is
    pure exact, and therefore
    \mbox{$\operatorname{l.\!p.\!gl.\!dim}R$} equals the (ordinary)
    left global dimension of $R$. This applies for example to the
    rings:
    \begin{itemlist}
    \item If $R$ has left global dimension zero, that is, $R$ is left
      semi-simple, then \mbox{$\operatorname{l.\!p.\!gl.\!dim}R=0$}.
    \item The ring \mbox{$R = \{ (x_n)_{n \in \mathbb{N}} \in
        k^\mathbb{N} |\, x_n \text{ constant for } n \gg 0\}$} ($k$
      any field), is von Neumann regular with unit.  Also, $R$ has
      global dimension $1$: As $R$ is not noetherian its global
      dimension is \mbox{$>\!0$}. Since $R$ is von Neumann regular,
      every ideal is generated by idempotents\footnote{This follows
        easily from the fact \cite[thm.\,p.\,10]{Goodearl} that every
        principal ideal is generated by an idempotent.}. Clearly, $R$
      has only $\aleph_0$ many idempotents, so
      \cite[cor.\,2.47]{Osofsky} implies the claim.
    \end{itemlist}
    For later use we note that a von Neumann regular ring is
    automatically coherent from either side.
  \end{rqm}
\end{bfhpg}

\begin{proof}[Proof of Example \eqref{exa:pure}]
  ``(a)'': In Definition \eqref{dfn:pure} we may take $d$ to be the
  left pure global dimension $d$ of $R$.

  \medskip
  \noindent
  ``(b)'': If every flat module has finite projective dimension then,
  in fact, there exists $d$ such that \mbox{$\mathrm{pd}_RF \leqslant
    d$} for all \mbox{$F \in \Flat\,R$}. If \mbox{$\mathbb{A}
    \subseteq \Proj\,R$} we get $\F(\mathbb{A}) \subseteq \Flat\,R$ by
  Lemma \eqref{lem:FProjandFInj}(a), and it follows immediately that
  the number $d$ implements the (\textsc{pp}) property for
  $\mathbb{A}$.
  
  \medskip
  \noindent
  ``(c)'': If \mbox{$\mathbb{A} \subseteq \Inj\,R$} we get
  \mbox{$\F(\mathbb{A}) \subseteq \Inj\,R$} by Lemma
  \eqref{lem:FProjandFInj}(b), and it follows that the number $d=0$
  implements the (\textsc{pi}) property for $\mathbb{A}$.
\end{proof}

\begin{proof}[Proof of Proposition \eqref{prp:PPandPI}]
  Let \mbox{$A = \cdots \to A_{n+1} \to A_n \to A_{n-1} \to \cdots$}
  be a pure exact sequence with $A_i \in \mathbb{A}$, and decompose
  $A$ into short exact sequences,
  \begin{displaymath}
    S_n = \ 0 \longrightarrow \Omega_n \longrightarrow A_n
    \longrightarrow \Omega_{n-1} \longrightarrow  0. 
  \end{displaymath}
  It follows that every $S_n$ is pure exact. We want to prove that
  $S_n$ is split exact, so it suffices to show that $\Omega_{n-1}$ is
  pure projective, or that $\Omega_n$ is pure injective. We will
  actually prove the following:
  \begin{prt}
  \item If $\mathbb{A}$ has property \textsc{(pp)} then every
    $\Omega_n$ is pure projective.
  \item If $\mathbb{A}$ has property \textsc{(pi)} then every
    $\Omega_n$ is pure injective.
  \end{prt}
  We will only prove (a), as the proof of (b) is similar: By
  Definition \eqref{dfn:F} every $\Omega_m$ belongs to
  $\F(\mathbb{A})$. To see that $\Omega_n$ is pure projective we
  consider the pure exact sequence,
  \begin{displaymath}
    0 \longrightarrow \Omega_n \longrightarrow A_n \longrightarrow
    \cdots \longrightarrow A_{n-d+1} \longrightarrow \Omega_{n-d}
    \longrightarrow 0,  
  \end{displaymath}
  where $d$ is a number which implements the property \textsc{(pp)}
  for $\mathbb{A}$. Since $\Omega_{n-d}$ belongs to $\F(\mathbb{A})$,
  and $A_{n-d+1}, \ldots, A_n$ belong to $\mathbb{A}$, the property
  \textsc{(pp)} guarantees that $\Omega_n$ is pure projective.
\end{proof}

\section{Summary} \label{sec:summary}

\noindent
Using the results from the previous sections we are now able to give a
list of examples of concrete module classes $\X=\mathsf{Add}\,\X$,
cf.~\eqref{bfhpg:setup}, such that the triangulated homotopy category
$\sK(\X)$ is compactly generated. In most of our examples, rings with
finite pure global dimension play an important role,
cf.~\eqref{bfhpg:pgldim}.  \smallskip

\begin{rqm}
\item Assume that $R$ is right coherent with finite left pure global
  dimension, and that $M$ is an $R$--bimodule which is finitely
  presented from either side. Then we may take
  \begin{displaymath}
    \X=\mathsf{Add}({}_RM).
  \end{displaymath}
  If \mbox{${}_RM$} is projective then \mbox{$\mathsf{Add}({}_RM)
    \subseteq \Proj\,R$}, and the condition ``finite left pure global
  dimension'' may be replaced by ``every flat module has finite
  projective dimension''. Thus, under this assumption the special case
  \mbox{$M=R$} recovers \cite[thm.\,2.4]{Jorgensen}, namely that we
  may take
  \begin{displaymath}
    \X=\Proj\,R.
  \end{displaymath}
  \emph{References:} \eqref{thm:main}, \eqref{prp:AddM},
  \eqref{prp:PPandPI}, and \eqref{exa:pure}(a)(b). \smallskip
  
\item Assume that $R$ is left noetherian with finite left pure global
  dimension, and that $M_1,\ldots,M_n$ are $R$--bimodules such that
  each $M_j$ is finitely generated as a left $R$--module. Then we may
  take 
  \begin{displaymath}
    \X=\mathsf{Inj}(M_1,\ldots,M_n).
  \end{displaymath}
  If every \mbox{$(M_j)_R$} is flat then
  \mbox{$\mathsf{Inj}(M_1,\ldots,M_n)\subseteq\Inj\,R$}, and the
  condition ``finite left pure global dimension'' is superfluous. In
  particular, the special case \mbox{$n=1$} and \mbox{$M_1=R$}
  recovers \cite[prop.\,2.3]{Krause}, namely that we may take
  \begin{displaymath}
    \X=\Inj\,R.
  \end{displaymath}
  \emph{References:} \eqref{thm:main}, \eqref{prp:Minj},
  \eqref{dfn:DefOfMinj}, \eqref{prp:PPandPI}, and
  \eqref{exa:pure}(a)(c).  \smallskip
  
\item If $R$ has finite left pure global dimension then we may take
  \begin{displaymath}
    \X=\PP\,R.
  \end{displaymath}
  In particular, if $R$ is left pure-semisimple then we
  can use 
  \begin{displaymath}
    \X=\Mod\,R.
  \end{displaymath}
  \emph{References:} \eqref{thm:main}, \eqref{exa:PROJ}(a),
  \eqref{exa:inj}(b), \eqref{prp:PPandPI}, and \eqref{exa:pure}(a).
  \smallskip
  
\item Assume that $R$ is commutative and noetherian with a duali\-zing
  complex, and that $R$ has finite pure global dimension. If
  \mbox{$\,\varinjlim\mathsf{Gproj}\,R=\mathsf{GFlat}\,R$} (this
  happens for example if $R$, in addition, has finite injective
  dimension over itself) then we may take
  \begin{displaymath}
    \X=\mathsf{Add}(\mathsf{Gproj}\,R).
  \end{displaymath}
  \emph{References:} \eqref{thm:main}, \eqref{prp:Gproj},
  \eqref{prp:PPandPI}, and \eqref{exa:pure}(a). \smallskip
  
\item If $R$ is left noetherian and right coherent with finite left
  pure global dimension then we may take for example
  \begin{align*}
    \X &= \mathsf{Add}(\mathsf{sub}(\Flat\,R)).
  \end{align*}
  \noindent
  \emph{References:} \eqref{thm:main}, \eqref{prp:Addsub},
  \eqref{bfhpg:apply}(b), \eqref{bfhpg:apply}(1), \eqref{prp:PPandPI},
  and \eqref{exa:pure}(a). \smallskip
  
\item If $R$ is commutative and noetherian with finite pure global
  dimension, and \mbox{$S\subseteq R$} is a multiplicative subset then
  we may take
  \begin{align*}
    \X 
      &= 
    \mathsf{Add}(\mathsf{sub} \{\text{$S$--torsion free
      modules}\}) \\ 
      &= 
    \mathsf{Add} \{\text{finitely generated $S$--torsion free
      modules}\}.
  \end{align*}
  \noindent
  \emph{References:} \eqref{thm:main}, \eqref{prp:Addsub},
  \eqref{bfhpg:apply}(c), \eqref{bfhpg:apply}(1), \eqref{prp:PPandPI},
  and \eqref{exa:pure}(a). \smallskip
  
\item Assume that $R$ is coherent from either side with finite left
  pure global dimension, and that $M$ is an $R$--bimodule which is
  finitely presented from either side.  Then we may take
  \begin{align*}
    \X &= \mathsf{Add}(\mathsf{sub}(\mathsf{add}({}_RM))).
  \end{align*}
  \noindent
  \emph{References:} \eqref{thm:main}, \eqref{prp:Addsub},
  \eqref{bfhpg:apply}(d), \eqref{bfhpg:apply}(2), \eqref{prp:PPandPI},
  and \eqref{exa:pure}(a).
\end{rqm}

\begin{prp} \label{prp:flat}
  Assume that $R$ is right coherent. Then $\Flat\,R$ is preenveloping,
  so it makes sense to consider the set of compact objects
  \begin{displaymath}
    \mathcal{G}_{\Flat\,R} \subseteq \sK(\Flat\,R)
  \end{displaymath}
  from Theorem \eqref{thm:main}. The set $\mathcal{G}_{\Flat\,R}$
  generates $\sK(\Flat\,R)$ if and only if $\Flat\,R = \Proj\,R$.
\end{prp}

\begin{proof}
  That $\Flat\,R$ is preenveloping over right coherent rings follows
  from \cite[prop.\,6.5.1]{EnochsRHA}.  If \mbox{$\Flat\,R =
    \Proj\,R$} then we know from Section \ref{sec:summary}(1) above
  that \mbox{$\sK(\Flat\,R)=\sK(\Proj\,R)$} is generated by
  \mbox{$\mathcal{G}_{\Flat\,R}=\mathcal{G}_{\Proj\,R}$}. When
  \mbox{$\Flat\,R \neq \Proj\,R$} there exists a flat module $F$ which
  is not projective. Let
  \begin{displaymath}
    \tag{\text{$*$}}
    \cdots \longrightarrow P_2 \longrightarrow P_1 \longrightarrow P_0
    \longrightarrow F \longrightarrow 0
  \end{displaymath}
  be an augmented projective resolution of $F$, and note that $(*)$ is
  pure exact but not split. Therefore Theorem \eqref{thm:main} implies
  that $\mathcal{G}_{\Flat\,R}$ does not generate $\sK(\Flat\,R)$.
\end{proof}


\medskip
\noindent
Of course, Proposition \eqref{prp:flat} above does not rule out the
possibility that $\sK(\Flat\,R)$ could be generated by some larger set
of compact objects than $\mathcal{G}_{\Flat\,R}$. Hence we pose the
following:

\begin{bfhpg}[Question]
  When is $\sK(\Flat\,R)$ compactly generated?
\end{bfhpg}

\section*{Acknowledgements}

\noindent
We sincerely thank Christian U. Jensen for his willingness to answer
questions about pure global dimension, and in particular for pointing
out \eqref{bfhpg:pgldim}(7).


\providecommand{\bysame}{\leavevmode\hbox to3em{\hrulefill}\thinspace}
\providecommand{\MR}{\relax\ifhmode\unskip\space\fi MR }
\providecommand{\MRhref}[2]{%
  \href{http://www.ams.org/mathscinet-getitem?mr=#1}{#2}
}
\providecommand{\href}[2]{#2}

\enlargethispage{15ex}

\end{document}